\documentclass[11pt, a4paper]{article}
\usepackage[left = 2cm, right = 2cm, top = 3.5cm, bottom = 2.5cm]{geometry}
\usepackage[utf8]{inputenc}
\DeclareUnicodeCharacter{FEFF}{}

\usepackage{fancyhdr}
\usepackage[numbers]{natbib}
\usepackage[]{hyperref}
\hypersetup{colorlinks = true}
\usepackage{amsmath, amssymb, amsthm}

\newcommand{\ds}{\\[13pt]}
\AtBeginDocument{%
  \setlength{\abovedisplayskip}{\dimexpr\abovedisplayskip}
  \setlength{\belowdisplayskip}{\dimexpr\belowdisplayskip}
  \setlength{\abovedisplayshortskip}{\dimexpr\abovedisplayshortskip}
  \setlength{\belowdisplayshortskip}{\dimexpr\belowdisplayshortskip}
}

\title{\Large A recurrence relation for the odd order moments of the Fabius function}
\author{Søren G. Have}

\begin{document} \noindent
\maketitle
\renewcommand{\abstractname}{\vspace{-\baselineskip}}
\begin{abstract} \noindent
A simple recurrence relation for the even order moments of the Fabius function is proven. Also, a very similar formula for the odd order moments in terms of the even order moments is proved. The matrices corresponding to these formulas (and their inverses) are multiplied so as to obtain a matrix that correspond to a recurrence relation for the odd order moments in terms of themselves. The theorem at the end gives a closed-form for the coefficients.
\end{abstract}
\hspace{1cm}\\
The Fabius function $F:\mathbb{R}_0^+\mapsto\mathbb{R}$ satisfies
\begin{align}
F(0)=0\quad,\quad F'(x)=2F(2x)\quad,\quad \forall x \in[0,1]\ :\ F(1+x)=1-F(1-x)
\label{basic}
\end{align}
We define $P_i:[0,1]\mapsto\mathbb{R}$ by 
\begin{align*}
P_i(x) = F(2^{-i}(1+x))-(-1)^i F(2^{-i}(1-x))
\end{align*}
The 3rd property in (\ref{basic}) implies $P_0(x)=0$. Moreover
\begin{align*}
P_i'(x)= & 2^{-i} F'(2^{-i} (1+x)) +(-1)^i 2^{-i} F'(2^{-i} (1-x)) =\\
& 2^{1-i}(F(2^{1-i}(1+x))-(-1)^{i-1}F(2^{1-i}(1-x)) = 2^{1-i}P_{i-1}(x)
\end{align*}
Hence
\begin{align}\label{integral}
P_{i+1}(x)=2^{-i} \int_0^x P_i(u)du +P_{i+1}(0)=
2^{-i} \int_0^x P_i(u)du +\begin{cases} 
      0 & i\text{ is odd} \\
      2F(2^{-i-1}) & \text{else} 
   \end{cases}
\end{align}
Let $c_i(j)$ denote the coefficient of $x^j$ in $P_i(x)$. It follows from (\ref{integral}) that
\begin{align*}
c_i(0) = \begin{cases} 
      0 & i\text{ is odd} \\
      2F(2^{-i-1}) & \text{else} 
   \end{cases} \quad,\quad c_i(j) = 
   \begin{cases} 
      2^{1-i}c(i-1,j-1)/j & j<i \\
      0 & \text{else} 
   \end{cases}
\end{align*}
Which has the solution
\begin{align*}
c_i(j) = \begin{cases} 
      \frac{2^{(1-i)j}F(2^{j-i})}{2^{(1-j)j/2-1} j!} & 0\leq j < i\text{ and } j-i \text{ is odd} \\
      0 & \text{else} 
   \end{cases}
\end{align*}
It holds that
\begin{align} \label{lem1}
\sum_{j=0}^{i+1} c_{i+2}(j) = P_{i+2}(1) = F(2^{-i-1})-(-1)^{i+2} F(0) =F(2^{-i-1})
\end{align}
Equating the RHS and LHS of (\ref{lem1}) we get for $i\in 2\mathbb{N}$ the following formula for $F(2^{-i-1})$:
\begin{align*}
F(2^{-i-1}) =& \sum_{j=0}^{i+1} c_{i+2}(j) = \sum_{j=0}^{i/2} c_{i+2}(2j+1) = \sum _{j=0}^{i/2}\frac{2^{(1 - (i + 2)) (2 j + 1)}F(2^{2 j + 1 - (i + 2)})}{2^{(1 - (2 j + 1)) (2 j + 1)/2 - 1}(2 j + 1)!}  =\\&
2^{-i}F(2^{-i-1})+ \sum _{j=1}^{i/2}\frac{4^{j(j-i)}F(2^{2 j - i - 1})}{2^{i+j}(2 j + 1)!}\ \Longrightarrow\ F(2^{-i-1}) = 
\sum _{j=1}^{i/2}\frac{4^{j(j-i)}F(2^{2 j - i - 1})}{(2^{i+j}-2^j)(2 j + 1)!}
\end{align*}
From the same equation we get for $i\in 2\mathbb{N}$ the following formula for $F(2^{-i})$:
\begin{align*}
F(2^{-i}) = \sum_{j=0}^i c_{i+1}(j) = \sum_{j=0}^{i/2} c_{i+1}(2 j) =
\sum_{j=0}^{i/2} \frac{2^{(1-(i+1)) 2j} F(2^{2 j-(i+1)})}{
2^{(1-2 j)(2 j)/2-1} (2 j)!} = 
\sum _{j=0}^{i/2} \frac{4^{j (j-i)} F\left(2^{2 j-i-1}\right)}{2^{j-1} (2 j)!}
\end{align*}
The formulas depends on $F(1/2)$ which equals $1/2$ cf. the 3rd property in (\ref{basic}).
\ds\ds
We define $d_i$ by
\begin{align*}
d_i=2^{i (i+1)/2}\hspace{0.4mm}i!\hspace{0.3mm}F(2^{-i-1});
\end{align*}
Clearly $d_0 = 1/2$, and if $i\in2\mathbb{N}-1$
\begin{align*}
d_i= & 2^{i(i+1)/2}\hspace{0.3mm} i! \sum_{j=0}^{(i+1)/2} \frac{4^{j(j-i-1)}F(2^{2j-i-2})}{2^{j-1}(2j)!} = \\ &
\sum_{j=0}^{(i+1)/2} \frac{2^{i (i+1)/2} i!}{2^{j-1}\hspace{0.3mm}(2j)!}\frac{4^{j(j-i-1)}d_{i-2j+1}}{2^{(i-2j+1)((i-2j+1)+1)/2}\hspace{0.3mm}(i-2j+1)!} =
\sum_{j=0}^{(i+1)/2}\binom{i+1}{2j} \frac{d_{2j}/2^i}{i+1}
\end{align*}
If $i\in2\mathbb{N}$ we get
\begin{align*}
d_i= & 2^{i(i+1)/2}\hspace{0.3mm} i! \sum_{j=1}^{i/2} \frac{4^{j(j-i)}F(2^{2j-i-1})}{(2^{i+j}-2^{j})(2j+1)!} = \\ &
\sum_{j=1}^{i/2} \frac{2^{i (i+1)/2} i!}{(2^{i+j}-2^{j})(2j+1)!}\frac{4^{j(j-i)}d_{i-2j}}{2^{(i-2j)((i-2j)+1)/2}\hspace{0.3mm}(i-2j)!} =
\sum_{j=0}^{i/2 - 1}\binom{i}{2j} \frac{d_{2j}/(2^i-1)}{i-2j+1}
\end{align*}
The following holds c.f. theorem 6 of \cite{k1}.
\begin{align} \label{mu}
\mu_n=\int_0^1 F(x)x^n dx=\frac{1}{n+1}-d_n
\end{align}
Hence the substitution $d_i = 1/(i+1) - \mu_i$ turns any recurrence relation for $d_i$ into a recurrence relation for the moments of the Fabius function.
\ds
For all $i\in\mathbb{N}_0$ the recurrence relations above imply the following matrix forms, where $M_i$ is an $i\times(i+1)$ matrix, $R_i$ is a length $i$ row, $e_1^i$ is the length $i$ unit vector $(1,0,...,0)$, and $I_i$ is the $i\times i$ identity matrix
\begin{align*} &
\left(
\begin{array}{c}
 d_1 \\
 d_3 \\
 \vdots \\
 d_{2 i-1}
\end{array}
\right) = M_i \left(
\begin{array}{c}
 d_0 \\
 d_2 \\
 \vdots \\
 d_{2 i}
\end{array}
\right)\ \Longrightarrow\ 
\left(
\begin{array}{c}
 2d_0 \\
 d_1 \\
 d_3 \\
 \vdots \\
 d_{2 i-1}
\end{array}
\right) = \left(
\begin{array}{c}
 2 e_1^{i+1} \\
 M_i \\
\end{array}
\right) \left(
\begin{array}{c}
 d_0 \\
 d_2 \\
 \vdots \\
 d_{2 i}
\end{array}
\right) \hspace{4cm} \\ &
\left(
\begin{array}{c}
 d_0 \\
 d_2 \\
 \vdots \\
 d_{2 i+2}
\end{array}
\right) =\left(
\begin{array}{c}
 I_{i+1} \\
 R_{i+1} \\
\end{array}
\right) \left(
\begin{array}{c}
 2 e_1^{i+1} \\
 M_i \\
\end{array}
\right)^{-1} \left(
\begin{array}{c}
 2d_0 \\
 d_1 \\
 d_3 \\
 \vdots \\
 d_{2 i-1}
\end{array}
\right)\ \Longrightarrow \hspace{3cm}\\ & \hspace{3cm}
\left(
\begin{array}{c}
 d_1 \\
 d_3 \\
 \vdots \\
 d_{2 i+1}
\end{array}
\right) =M_{i+1}\left(
\begin{array}{c}
 I_{i+1} \\
 R_{i+1} \\
\end{array}
\right) \left(
\begin{array}{c}
 2 e_1^{i+1} \\
 M_i \\
\end{array}
\right)^{-1} \left(
\begin{array}{c}
 2d_0 \\
 d_1 \\
 d_3 \\
 \vdots \\
 d_{2 i-1}
\end{array}
\right)\
\end{align*}
Invertibility is no issue, because a pseudo-inverse suffices. The last matrix form expresses $d_{2i+1}$ as a linear combination of previous odd-indexed $d_i$s and $2d_0=1$. The coefficients are the last row of the matrix:
\begin{align*}
G_{i+1}=M_{i+1}\left(
\begin{array}{c}
 I_{i+1} \\
 R_{i+1} \\
\end{array}
\right) \left(
\begin{array}{c}
 2 e_1^{i+1} \\
 M_i \\
\end{array}
\right)^{-1}
\end{align*}
The $(k,j)$'th entry of $M_i$, and the $j$'th entry of $R_i$ are known from the recurrence relations:
\begin{align*}
(M_i)_{k,j} = \binom{2 k}{2 j-2}\frac{1}{4^k k} \quad,\quad
(R_i)_j = \binom{2i}{2j-2}\frac{1/(4^i-1)}{2i-2j+3}
\end{align*}
\paragraph{Theorem:} For all $i\in\mathbb{N}$ the $j$'th column entry of the last row of $G_i$ is
\begin{align*}
(G_i)_{i,j} = \left(\frac{E_{2 (i-j+1)}}{4^{2-j}} + 
\frac{\zeta (2 j-2 i-3)}{2 j-2 i-3} (1-4^{i-j+2})\right)\frac{4}{1-4^i}\begin{cases} 
      \binom{2i-1}{2j-3} & j\neq1 \\
      1/(2i) & \text{else} 
   \end{cases}
\end{align*}
And so for all $i\in2\mathbb{N}-1$ the following recurrence relation holds:
\begin{align*}
d_i =  \sum_{\substack{j\ \in\ 2\mathbb{N}_0 - 2 \\ j\ <\ i-1\ \ \ }}
 \left(\frac{\zeta (j-i)}{j-i}-\frac{2^j E_{i-j-1}}{2^{i-j+1}-1}\right)\frac{2^{i-j+3}-4}{2^{i+1}-1}
 \begin{cases} 
      \binom{i}{j+1}d_{j+1} & j\neq -2 \\
      1/(i + 1) & \text{else} 
   \end{cases}
\end{align*}
\paragraph{Proof:} It suffices to show that the last row of $G_i \left(
\begin{array}{c}
 2 e_1^{i} \\
 M_{i-1} \\
\end{array}
\right)$ and $M_{i} \left(
\begin{array}{c}
 I_{i} \\
 R_{i} \\
\end{array}
\right)$ are the same for $i\in\mathbb{N}$.\ds
Each entry in the last row for the second matrix product consists of just 2 terms:
\begin{align*}
\left(M_i \left(
\begin{array}{c}
 I_i \\
 R_i \\
\end{array}
\right)\right)_{i,j} = (M_i)_{i,j} + (M_i)_{i,i+1}(R_i)_j =\ &
\binom{2 i}{2 j-2}\frac{1}{4^i i} + 
\frac{1}{4^i i} \binom{2 i}{2j-2} \frac{1/(4^i-1)}{2 i - 2 j + 3} = \\ &
\binom{2 i-1}{2 j-2} \frac{2^{1-2 i}(j-i-1)+2 i-2 j+3}{
(1-4^i)(j-i-1)(2 i-2 j+3)}
\end{align*}
The last row for the first matrix product:
\begin{align*}
& \left(G_i \left(
\begin{array}{c}
 2 e_1^i \\
 M_{i-1} \\
\end{array}
\right)\right)_{i,j} = \left(\sum_{m=2}^i (G_i)_{i,m} (M_{i-1})_{m-1,j}\right) + \begin{cases} 
      2 (G_i)_{i,1} & j=1 \\
      0 & \text{else} 
   \end{cases} = \\ &
\sum_{m=2}^i \left(\frac{E_{2 (i-m+1)}}{4^{2-m}} + 
\frac{\zeta (2m-2i-3)}{2m-2i-3} (1-4^{i-m+2})\right)
\frac{4}{1-4^i} \binom{2 i-1}{2 m-3}
\binom{2 m-2}{2 j-2}\frac{4^{1-m}}{m-1}\ +\\ & \hspace{0.4cm}
\left(\frac{E_{2 i}}{4}+\frac{\zeta(-2 i-1)}{-2 i-1} (1-4^{i+1})\right)\frac{4}{1-4^i}\frac{1}{i} \begin{cases} 
      1 & j=1 \\
      0 & \text{else} 
   \end{cases}
\end{align*}
At this point the last row for both matrix products are divided by $Q=\frac{4}{1-4^i}\binom{2 i-1}{2 j-2}$:
\begin{align*}
\frac{1}{Q}\left(M_i \left(
\begin{array}{c}
 i_i \\
 R_i \\
\end{array}
\right)\right)_{i,j} =\frac{2^{1-2 i} (j-i-1)+2 i-2 j+3}{
4(j-i-1)(2 i-2 j+3)}
\end{align*}
\begin{align*}
& \frac{1}{Q}  \left(G_i \left(
\begin{array}{c}
 2 e_1^i \\
 M_{i-1} \\
\end{array}
\right)\right)_{i,j}  = 
\left(\frac{E_{2 i}}{4}+\frac{\zeta (-2 i-1)}{-2 i-1}(1-4^{i+1})\right) \frac{1}{i}\begin{cases} 
      1 & j=1 \\
      0 & \text{else} 
   \end{cases}\ +\ \\& \hspace{0.2cm}
\sum_{m=2}^i \left(\frac{E_{2 (i-m+1)}}{4^{2-m}} + 
\frac{\zeta (2m-2i-3)}{2m-2i-3} (1-4^{i-m+2})\right)
\binom{2 i-1}{2 m-3}
\binom{2 m-2}{2 j-2}\frac{4^{1-m}}{m-1}/
\binom{2 i-1}{2 j-2} = \\& \hspace{4.4cm}
 \sum_{m=1}^i \left(\frac{E_{2 (i-m+1)}}{4^{2-m}} + 
\frac{\zeta (2m-2i-3)}{2m-2i-3} (1-4^{i-m+2})\right)
\frac{4^{1-m}}{i-m+1}\binom{2 i-2 j+1}{2 m-2 j}
\end{align*}
Next $j$ is substituted by $i-j+1$ and $m$ is substituted by $i-m+1$:
\begin{align*}
& \frac{1}{Q}\left(M_i \left(
\begin{array}{c}
 i_i \\
 R_i \\
\end{array}
\right)\right)_{i,j} = \frac{2^{-2 i-1}}{2 j+1}-\frac{1}{4 j} \\ &
\frac{1}{Q}\left(G_i \left(
\begin{array}{c}
 2 e_1^i \\
 M_{i-1} \\
\end{array}
\right)\right)_{i,j} =
\sum _{m=1}^i \left(\frac{E_{2 m}}{4^{m-i+1}}+\frac{\zeta (-1-2 m)}{-1-2 m}(1-4^{m+1})\right)\frac{4^{m-i}}{m}\binom{2 j-1}{2 j-2 m}
\end{align*}
The rows are subtracted from one another and the $\zeta$-function is expressed with Bernoulli numbers:
\begin{align*}
& \frac{1}{4 j} - \frac{2^{-2 i-1}}{2 j+1} +
\sum _{m=1}^i \left(\frac{E_{2 m}}{4^{m-i+1}}+\frac{B_{2 (m+1)}}{2 (m+1) (2 m+1)}(1-4^{m+1})\right)\frac{4^{m-i}}{m}\binom{2 j-1}{2 j-2 m} = \\ &
\left(\frac{1}{4 j}+\sum _{m=1}^i \frac{E_{2 m}}{4 m }\binom{2 j-1}{2 j-2 m}\right) - \left(\frac{2^{-2 i-1}}{2 j+1}-\sum_{m=1}^i \frac{B_{2 (m+1)} (4^{m-i}-4^{2 m-i+1})}{2 m (m+1) (2 m+1)}\binom{2 j-1}{2 j-2 m}\right)
\end{align*}
We show that both parentheses equal 0. Simplifying the first parenthesis yields the summation:
\begin{align*}
\frac{1}{4 j} + \sum_{m=1}^i \frac{E_{2 m}}{(4 m)}\binom{2 j}{2 m}\frac{m}{j} = 1 + \sum_{m=1}^i E_{2m}\binom{2 j}{2 m} = \sum_{m=0}^j E_{2m}\binom{2 j}{2 m}
\end{align*}
The sum is rewritten as the coefficient of a generating function:
\begin{align*}
\sum_{j=0}^{\infty} & \frac{z^{2j}}{(2j)!}\sum_{m=0}^j
\frac{E_{2m}(2j)!}{(2m)!(2j-2m)!} = \\ &
\sum_{j=0}^{\infty} \sum_{m=0}^j \frac{E_{2m} z^{2j}}{(2m)!(2j-2m)!} = \sum_{m=0}^{\infty} \sum_{j=m}^{\infty} \frac{E_{2m} z^{2j}}{(2m)!(2j-2m)!} = \\ &
\sum_{m=0}^{\infty} \frac{E_{2m}z^{2m}}{(2m)!} \sum_{j=0}^{\infty} \frac{z^{2j}}{(2j)!} = \cosh(z) \text{sech}(z) = 1 =
\sum_{j=0}^{\infty} \frac{z^{2j}}{(2j)!} \begin{cases} 
      1 & j=0 \\
      0 & \text{else} 
   \end{cases}
\end{align*}
Hence from linear independence of the monomials we have that the first parenthesis is 0 for $j\in\mathbb{N}$.\ds
Rewriting the second parenthesis makes it
\begin{align*}
1 + \sum _{m=0}^j \frac{B_{2 m+2}(4^{m+1}-16^{m+1}) (2 j+1)!}{
(2 j-2 m)!(2 m+2)!}
\end{align*}
the sum term of which is used as the coefficient of a generating function:
\begin{align*}
& \sum_{j=0}^{\infty} \frac{z^{2 j+1}}{(2 j+1)!} \sum_{m=0}^{j}\frac{B_{2 m+2}(4^{m+1}-16^{m+1}) (2 j+1)!}{
(2 j-2 m)!(2 m+2)!} =
\sum_{j=0}^{\infty}\sum_{m=0}^{j} \frac{z^{2 j+1}B_{2 m+2}(4^{m+1}-16^{m+1})}{
(2 j-2 m)!(2 m+2)!}= \\
& \sum_{m=0}^{\infty} \frac{B_{2 m+2}(4^{m+1}-16^{m+1})}{
(2 m+2)!} \sum_{j=m}^{\infty} \frac{z^{2 j+1}}{(2 j-2 m)!} = 
\sum_{m=0}^{\infty} \frac{B_{2 m+2}(4^{m+1}-16^{m+1})z^{2m+1}}{
(2 m+2)!} \sum_{j=0}^{\infty} \frac{z^{2 j}}{(2 j)!} = \\ &
\cosh(z)\frac{1}{2 z} \sum_{m=0}^{\infty}\left(
\frac{B_m (2 z)^m}{m!} + \frac{B_m (-2 z)^m}{m!} - 
\frac{B_m (4 z)^m}{m!} - \frac{B_m (-4 z)^m}{m!}
\right) = \\ &
\cosh(z)\frac{1}{2 z}\left(
\frac{z}{e^{2 z}-1} -\frac{z}{e^{-2 z}-1} -
\frac{2z}{e^{4 z}-1} +\frac{2z}{e^{-4 z}-1}
\right) = -\sinh (z)= \sum _{j=0}^{\infty } \frac{z^{2 j+1}}{(2 j+1)!}(-1)
\end{align*}
Hence from linear independence of the monomials we have that
\begin{align*}
1 + \sum _{m=0}^j \frac{B_{2 m+2}(4^{m+1}-16^{m+1}) (2 j+1)!}{
(2 j-2 m)!(2 m+2)!} = 1-1=0
\end{align*}
\qed

\end{document}